\documentclass[12pt]{article}%
\usepackage{amsmath}%
\usepackage{amsfonts}%
\usepackage{amssymb}%
\usepackage{graphicx}
\providecommand{\U}[1]{\protect\rule{.1in}{.1in}}
\newtheorem{theorem}{Theorem}[section]

\newtheorem{definition}[theorem]{Definition}

\newtheorem{problem}[theorem]{Problem}

\numberwithin{equation}{section}

\begin{document}

\title{Reconstruction techniques for quantum trees}
\author{Sergei A. Avdonin{\small $^{\text{1}}$}, Kira V. Khmelnytskaya$^{\text{2}}$,
Vladislav V. Kravchenko{\small $^{\text{3}}$}\\{\small $^{\text{1}}$ Department of Mathematics and Statistics, University of
Alaska, Fairbanks, AK 99775, USA}\\{\small $^{\text{2}}$ Faculty of Engineering, Autonomous University of
Queretaro, }\\{\small Cerro de las Campanas s/n, col. Las Campanas Quer\'{e}taro, Qro. C.P.
76010 M\'{e}xico}\\{\small $^{\text{3}}$ Department of Mathematics, Cinvestav, Campus
Quer\'{e}taro, }\\{\small Libramiento Norponiente \#2000, Fracc. Real de Juriquilla,
Quer\'{e}taro, Qro., 76230 M\'{e}xico}\\{\small e-mail: s.avdonin@alaska.edu, khmel@uaq.edu.mx,
vkravchenko@math.cinvestav.edu.mx}}
\maketitle

\begin{abstract}
The inverse problem of recovery of a potential on a quantum tree graph from
Weyl's matrix given at a number of points is considered. A method for its
numerical solution is proposed. The overall approach is based on the leaf
peeling method combined with Neumann series of Bessel functions (NSBF)
representations for solutions of Sturm-Liouville equations. In each step, the
solution of the arising inverse problems reduces to dealing with the NSBF coefficients.

The leaf peeling method allows one to localize the general inverse problem to
local problems on sheaves, while the approach based on the NSBF
representations leads to splitting the local problems into two-spectra inverse
problems on separate edges and reduce them to systems of linear algebraic
equations for the NSBF coefficients. Moreover, the potential on each edge is
recovered from the very first NSBF coefficient. The proposed method leads to
an efficient numerical algorithm that is illustrated by numerical tests.

\end{abstract}

\section{Introduction}

Quantum graphs and in particular inverse problems on\ quantum graphs find
numerous applications in science and engineering and give rise to challenging
problems involving many areas of modern mathematics, from combinatorics to
partial differential equations and spectral theory. A number of surveys and
collections of papers on quantum graphs appeared last years, including the
first books on this topic by Berkolaiko and Kuchment \cite{BerkolaikoKuchment}
and Mugnolo \cite{Mugnolo}. Theory of inverse problems on quantum graphs is an
actively developing research area of applied mathematics, providing existence
and uniqueness results and possible approaches for solution (see, e.g.,
\cite{Kurasov et al 2005, Yurko2005, Belishev Vakulenko 2006,
AvdoninKurasov2008, AvdoninLeugeringMikhaylov2010, AvdoninBell2015}).

To date, there are few papers presenting methods apt for numerical solution of
inverse problems on quantum graphs \cite{Belishev Vakulenko 2006,
AvdoninBelinskiyMatthews,AvdoninKravchenko2022, AKK2023}. Since the problems
of space discretization of differential equations on metric graphs turn out to
be very difficult, and even the forward boundary value problems on graphs
present considerable numerical challenges (see, e.g. \cite{ArioliBenzi2018}),
a promising path for efficient solution of inverse problems seems to go
through and involve analytical representations of solutions. In
\cite{AvdoninKravchenko2022, AKK2023} such approach was developed for star
shaped graphs. It is based on the use of very special functional series
representations for solutions of the Sturm-Liouville equation
\[
-y^{\prime\prime}+q(x)y=\rho^{2}y,\text{\quad}x\in\left(  0,L\right)  ,\quad
q\in L_{1}\left(  0,L\right)  ,\quad\rho\in\mathbb{C}.
\]
These representations have the form of so-called Neumann series of Bessel
functions (NSBF) (see, e.g., \cite[Chapter XVI]{Watson}, \cite{Wilkins} and
the recent monograph on the subject \cite{Baricz et al Book} and references
therein) and were obtained in \cite{KNT} as a result of the expansion of the
transmutation operator kernel (see, e.g., \cite{KarapetyantsKravchenkoBook},
\cite{KrBook2020}, \cite{LevitanInverse}, \cite{Marchenko},
\cite{SitnikShishkina Elsevier}) into a Fourier-Legendre series. The NSBF
representations possess certain unique features, which make them especially
convenient for solving inverse spectral problems. The remainders of the series
admit bounds, which are independent of $\operatorname{Re}\rho$. This
facilitates dealing with approximate solutions (partial sums of the NSBF) on
very large intervals in\ $\rho$ and with a non-deteriorating accuracy.
Moreover, the knowledge of the very first coefficient of the NSBF
representation allows one to recover the potential of the Sturm-Liouville
equation, that in practice ensures satisfactory numerical results even when a
reduced number of terms of the series is considered and consequently a reduced
number of linear algebraic equations is solved in each step.

In the present work we extend this approach, based on the NSBF
representations, onto an arbitrary quantum tree. This is done by applying the
leaf peeling method developed in \cite{AvdoninKurasov2008}, combined with new
ideas regarding the computation of solutions and their derivatives at
abscission points. The inverse problem considered here is the approximate
recovery of a potential on a quantum tree from a Weyl matrix given at a set of
points $\rho=\rho_{k}$, $k=1,\ldots,K$. The Weyl matrix of a quantum graph is
in fact its Dirichlet-to-Neumann map. It plays essential role in all aspects
of study of quantum graphs, including spectral theory \cite{Yurko2005} and
controllability \cite{AvdoninKurasov2008}. The physical meaning of the inverse
problem is the recovery of the differential operator from a measured response
at a number of frequencies of a system modelled by the quantum graph.

The leaf peeling method allows one to localize the inverse problem on a tree
to an inverse problem on a sheaf, which is a star shaped subgraph whose all
but one of the edges are leaf (boundary) edges. After the potential on the
leaf edges of the sheaf is recovered, they can be removed, and the leaf
peeling method allows one to calculate the Weyl matrix for the new smaller
tree. This procedure repeated finitely many times gradually exhausts the whole
tree, leading in the last step to an inverse problem on a star shaped graph. A
method for its solution with the aid of the NSBF representations was proposed
in \cite{AKK2023}, so that in the present paper we develop a method for
solving the local inverse problems on sheaves as well as techniques for
computing additional auxiliary functions required for applying the leaf
peeling method.

In fact, the overall approach consists in reducing the inverse problem on a
graph to operations with the NSBF coefficients of solutions on edges. When
solving the local problem, in the first step, we compute the NSBF coefficients
for a couple of linearly independent solutions on the leaf edges at the end
point of each edge, which is associated with the common vertex of the sheaf.
This first step allows us to split the problem into separate problems on leaf
edges. Second, the obtained NSBF coefficients are used for computing the
Dirichlet-Dirichlet and Neumann-Dirichlet eigenvalues of the potential on each
leaf edge. The first feature of the NSBF representations mentioned above
ensures the possibility to compute hundreds of the eigenvalues with a uniform
accuracy. Thus, on each leaf edge of the sheaf we obtain a two spectra inverse
Sturm-Liouville problem. Results on the uniqueness and solvability of such
problems are well known and can be found, e.g., in \cite{Chadan et al 1997},
\cite{LevitanInverse}, \cite{SavchukShkalikov}, \cite{Yurko2007}. To this
problem we apply the method from \cite{Kr2022Completion} which again involves
the NSBF representations. It provides for computing certain multiplier
constants \cite{Brown et al 2003} which relate to the Neumann-Dirichlet
eigenfunctions associated to the same eigenvalues but normalized at the
opposite endpoints of the interval. This leads to a system of linear algebraic
equations for the coefficients of their NSBF representations, already for
interior points of the interval. Solving the system we find the very first
coefficient, from which the potential is recovered. It is worth mentioning
that the NSBF representations were first used for solving inverse
Sturm-Liouville problems on a finite interval in \cite{Kr2019JIIP}. Later on,
the approach from \cite{Kr2019JIIP} was improved in \cite{KrBook2020} and
\cite{KT2021 IP1}, \cite{KT2021 IP2}. In those papers the system of linear
algebraic equations was obtained with the aid of the Gelfand-Levitan integral
equation. In \cite{Kr2022Completion} another approach, based on the
consideration of the eigenfunctions normalized at the opposite endpoints, was
developed, and this idea was used in \cite{AvdoninKravchenko2022},
\cite{KKC2022Mathematics}, \cite{AKK2023} and is used in the present work when
solving the two-spectra inverse problems on the edges.

Thus, the main result of this work is an efficient method for solving the
inverse problem on a quantum tree, consisting in the recovery of a potential
from a Weyl matrix. We discuss the numerical implementation of the method and
give numerical examples.

In Section \ref{Sect ProblemSetting} we recall the definition of the Weyl
matrix and formulate the inverse problem. In Section \ref{Sect FundSys} we
write the Weyl solutions in terms of the fundamental systems of solutions on
each edge and recall the NSBF representations for solutions of the
Sturm-Liouville equation as well as some of their relevant features. In
Section \ref{Sect LeafPeeling} we explain in detail the leaf peeling method.
In Section \ref{Sect Local Problem} we give a detailed description of the
proposed method for the solution of the local inverse problem. In Section
\ref{Sect SummaryMethod} we summarize\ the overall method for solving the
inverse problem on a tree. In Section \ref{Sect Numerical} we discuss the
numerical implementation of the method and numerical examples. Finally,
Section \ref{Sect Concl} contains some concluding remarks.

\section{Problem setting\label{Sect ProblemSetting}}

Let $\Omega$ be a finite connected compact graph without cycles (a tree graph)
consisting of $P$ edges $e_{1}$,...,$e_{P}$ and $P+1$ vertices $V=\left\{
v_{1},...,v_{P+1}\right\}  $. The notation $e_{j}\sim v$ means that the edge
$e_{j}$ is incident to the vertex $v$. Every edge $e_{j}$ is identified with
an interval $(0,L_{j})$ of the real line. The boundary $\Gamma=\left\{
\gamma_{1},\ldots,\gamma_{m}\right\}  $ of $\Omega$ is the set of all leaves
of the graph (the external vertices). The edge adjacent to some $\gamma_{j}%
\ $is called a leaf or boundary edge.

Let $q\in L_{2}(\Omega)$ be real valued, and $\lambda$ a complex number. A
continuous function $u$ defined on the graph $\Omega$ is a $P$-tuple of
functions $u_{j}\in C\left[  0,L_{j}\right]  $ satisfying the continuity
condition at the internal vertices $v$: $u_{i}(v)=u_{j}(v)$ for all
$e_{i},e_{j}\sim v$. Then $u\in C(\Omega)$.

We say that a function $u$ is a solution of the equation
\begin{equation}
-u^{\prime\prime}(x)+q(x)u(x)=\lambda u(x) \label{Schr}%
\end{equation}
on the graph $\Omega$ if besides (\ref{Schr}) the following conditions are
satisfied
\begin{equation}
u\in C(\Omega) \label{continuity}%
\end{equation}
and
\begin{equation}
\sum_{e_{j}\sim v}\partial u_{j}(v)=0,\text{\quad for all }v\in V\setminus
\Gamma. \label{KN}%
\end{equation}
Here $\partial u_{j}(v)$ stands for the derivative of $u$ at the vertex $v$
taken along the edge $e_{j}$ in the direction outward the vertex. The sum in
(\ref{KN}) is taken over all the edges incident to the internal vertex $v$.
Condition (\ref{KN}) is known as the Kirchhoff-Neumann condition.

Consider a solution $w_{i}$ of (\ref{Schr}) such that $w_{i}$ equals zero at
all leaves but one: $\gamma_{i}$, at which it equals one. That is $w_{i}$ is a
solution of (\ref{Schr}) such that
\[
w_{i}(\gamma_{i})=1\quad\text{and}\quad w_{i}(\gamma_{j})=0\text{ for all
}j\neq i.
\]
Such a solution $w_{i}$ is called the \textbf{Weyl solution} associated with
the leaf $\gamma_{i}$. Since the Dirichlet spectrum of (\ref{Schr}) is real,
the Weyl solution for any $\gamma_{i}$ exists and is unique for all
$\lambda\notin\mathbb{R}$.

\begin{definition}
The $m\times m$ matrix-function $\mathbf{M}(\lambda)$, $\lambda\notin
\mathbb{R}$, consisting of the elements $\mathbf{M}_{ij}(\lambda)=\partial
w_{i}(\gamma_{j})$, $i,j=1,\ldots,m$ is called the \textbf{Weyl matrix}.
\end{definition}

In fact, for a fixed value of $\lambda$, the Weyl matrix represents a
Dirichlet-to-Neumann map of the quantum graph defined by $\Omega$ and $q\in
L_{2}(\Omega)$. Indeed, if $u$ is a solution of (\ref{Schr}) satisfying the
Dirichlet condition at the boundary vertices $u(\gamma,\lambda)=f(\lambda)$,
then $\partial u(\gamma,\lambda)=\mathbf{M}(\lambda)f(\lambda)$,
$\lambda\notin\mathbb{R}$.

The problem we consider in the present paper can be formulated as follows.

\begin{problem}
Given $\Omega$ and the Weyl matrix at a finite number of points $\lambda_{k}$,
$k=1,\ldots,K$, find the potential $q(x)$ approximately.
\end{problem}

When the Weyl matrix or even its main diagonal is known everywhere, the
potential $q(x)$ is determined uniquely (see, e.g., \cite{Yurko2005}). The
knowledge of the Weyl matrix at a finite number of points may allow one to
recover the potential $q(x)$ only approximately.

Practical importance of this inverse problem is quite obvious. The entries of
the Weyl matrix represent the response of the physical system described by the
quantum graph to a unitary impulse applied at one end while isolating the
others. From the knowledge of this response at some values $\lambda_{k}$,
$k=1,\ldots,K$ we recover approximately the Sturm-Liouville equation
(\ref{Schr}) on the whole tree graph.

\section{Fundamental system of solutions and the Weyl
solutions\label{Sect FundSys}}

By $\varphi_{i}(\rho,x)$ and $S_{i}(\rho,x)$ we denote the solutions of the
equation
\begin{equation}
-y^{\prime\prime}(x)+q_{i}(x)y(x)=\rho^{2}y(x),\quad x\in(0,L_{i})
\label{Schri}%
\end{equation}
satisfying the initial conditions
\[
\varphi_{i}(\rho,0)=1,\quad\varphi_{i}^{\prime}(\rho,0)=0,
\]%
\[
S_{i}(\rho,0)=0,\quad S_{i}^{\prime}(\rho,0)=1.
\]
Here $q_{i}(x)$ is the component of the potential $q(x)$ on the edge $e_{i}$,
and $\rho=\sqrt{\lambda}$, $\operatorname{Im}\rho\geq0$. For a leaf edge
$e_{i}$ it is convenient to identify its leaf $\gamma_{i}$ with zero. Then the
Weyl solution $w_{i}(\rho,x)$ has the form%
\[
w_{i}(\rho,x)=\varphi_{i}(\rho,x)+\mathbf{M}_{ii}(\rho^{2})S_{i}(\rho
,x)\quad\text{on the adjacent edge }e_{i}%
\]
and
\[
w_{i}(\rho,x)=\mathbf{M}_{ij}(\rho^{2})S_{j}(\rho,x)\quad\text{on every leaf
edge }e_{j},\quad j\neq i.
\]
On internal edges $e_{j}$ we have
\[
w_{i}(\rho,x)=a_{ij}(\rho)\varphi_{j}(\rho,x)+b_{ij}(\rho)S_{j}(\rho,x),
\]
where the choice of which vertex lies at zero is arbitrary, and in general the
factors $a_{ij}(\rho)$, $b_{ij}(\rho)$ are unknown.

\begin{theorem}
[\cite{KNT}]\label{Th NSBF} The solutions $\varphi_{i}(\rho,x)$ and
$S_{i}(\rho,x)$ of (\ref{Schri}) and their derivatives with respect to $x$
admit the following series representations
\begin{align}
\varphi_{i}(\rho,x)  &  =\cos\left(  \rho x\right)  +\sum_{n=0}^{\infty
}(-1)^{n}g_{i,n}(x)\mathbf{j}_{2n}(\rho x),\label{phiNSBF}\\
S_{i}(\rho,x)  &  =\frac{\sin\left(  \rho x\right)  }{\rho}+\frac{1}{\rho}%
\sum_{n=0}^{\infty}(-1)^{n}s_{i,n}(x)\mathbf{j}_{2n+1}(\rho x),\label{S}\\
\varphi_{i}^{\prime}(\rho,x)  &  =-\rho\sin\left(  \rho x\right)  +\frac
{\cos\left(  \rho x\right)  }{2}\int_{0}^{x}q_{i}(t)\,dt+\sum_{n=0}^{\infty
}(-1)^{n}\gamma_{i,n}(x)\mathbf{j}_{2n}(\rho x),\label{phiprimeNSBF}\\
S_{i}^{\prime}(\rho,x)  &  =\cos\left(  \rho x\right)  +\frac{\sin\left(  \rho
x\right)  }{2\rho}\int_{0}^{x}q_{i}(t)\,dt+\frac{1}{\rho}\sum_{n=0}^{\infty
}(-1)^{n}\sigma_{i,n}(x)\mathbf{j}_{2n+1}(\rho x), \label{Sprime}%
\end{align}
where $\mathbf{j}_{k}(z)$ stands for the spherical Bessel function of order
$k$, $\mathbf{j}_{k}(z):=\sqrt{\frac{\pi}{2z}}J_{k+\frac{1}{2}}(z)$ (see,
e.g., \cite{AbramowitzStegunSpF}). The coefficients $g_{i,n}(x)$, $s_{i,n}%
(x)$, $\gamma_{i,n}(x)$ and $\sigma_{i,n}(x)$ can be calculated following a
simple recurrent integration procedure (see \cite{KNT} or \cite[Sect.
9.4]{KrBook2020}), starting with
\begin{align}
g_{i,0}(x)  &  =\varphi_{i}(0,x)-1,\quad s_{i,0}(x)=3\left(  \frac{S_{i}%
(0,x)}{x}-1\right)  ,\label{beta0}\\
\gamma_{i,0}(x)  &  =g_{i,0}^{\prime}(x)-\frac{1}{2}\int_{0}^{x}%
q_{i}(t)\,dt,\quad\sigma_{i,0}(x)=\frac{s_{i,0}(x)}{x}+s_{i,0}^{\prime
}(x)-\frac{3}{2}\int_{0}^{x}q_{i}(t)\,dt.\nonumber
\end{align}
For every $\rho\in\mathbb{C}$ all the series converge pointwise. For every
$x\in\left[  0,L_{i}\right]  $ the series converge uniformly on any compact
set of the complex plane of the variable $\rho$, and the remainders of their
partial sums admit estimates independent of $\operatorname{Re}\rho$. In
particular, for the partial sums
\[
\varphi_{i,N}(\rho,x):=\cos\left(  \rho x\right)  +\sum_{n=0}^{N}%
(-1)^{n}g_{i,n}(x)\mathbf{j}_{2n}(\rho x)
\]
and
\[
S_{i,N}(\rho,x):=\frac{\sin\left(  \rho x\right)  }{\rho}+\frac{1}{\rho}%
\sum_{n=0}^{N}(-1)^{n}s_{i,n}(x)\mathbf{j}_{2n+1}(\rho x)
\]
we have the estimates
\begin{equation}
\left\vert \varphi_{i}(\rho,x)-\varphi_{i,N}(\rho,x)\right\vert \leq
\frac{2\varepsilon_{i,N}(x)\,\sinh(Cx)}{C}\quad\text{and}\quad\left\vert
S_{i}(\rho,x)-S_{i,N}(\rho,x)\right\vert \leq\frac{2\varepsilon_{i,N}%
(x)\,\sinh(Cx)}{C} \label{estim S}%
\end{equation}
for any $\rho\in\mathbb{C}$ belonging to the strip $\left\vert
\operatorname{Im}\rho\right\vert \leq C$, $C\geq0$, where $\varepsilon
_{i,N}(x)$ is a positive function tending to zero as $N\rightarrow\infty$.
Analogous estimates are valid for the derivatives.
\end{theorem}

Roughly speaking, the approximate solutions and their derivatives approximate
the exact ones equally well for small and for large values of
$\operatorname{Re}\rho$. This is especially convenient when considering direct
and inverse spectral problems that requires operating on a large range of the
parameter $\rho$. This unique feature of the series representations
(\ref{phiNSBF})-(\ref{Sprime}) is due to the fact that they originate from an
exact Fourier-Legendre series representation of the integral kernel of the
transmutation operator \cite{KNT}, \cite[Sect. 9.4]{KrBook2020} (for the
theory of transmutation operators we refer to \cite{LevitanInverse},
\cite{Marchenko}, \cite{SitnikShishkina Elsevier}, \cite{Yurko2007}).

Moreover, the following statement is valid.

\begin{theorem}
\label{Th closeness of zeros}For any $\varepsilon>0$ there exists such
$N\in\mathbb{N}$ that all zeros of the functions $S_{i}(\rho,L_{i})$ and
$\varphi_{i}(\rho,L_{i})$ are approximated by corresponding zeros of the
functions $S_{i,N}(\rho,L_{i})$ and $\varphi_{i,N}(\rho,L_{i})$, respectively,
with errors uniformly bounded by $\varepsilon$. Moreover, $S_{i,N}(\rho
,L_{i})$ and $\varphi_{i,N}(\rho,L_{i})$, have no other zeros.
\end{theorem}

The proof of this statement is completely analogous to the proof of
Proposition 7.1 in \cite{KT2015JCAM} and consists in using (\ref{estim S}),
properties of characteristic functions of regular Sturm-Liouville problems and
the Rouch\'{e} theorem.

Series of the type
\[
\sum_{n=0}^{\infty}a_{n}J_{\nu+n}(z)
\]
are called Neumann series of Bessel functions (NSBF).

Note that formulas (\ref{beta0}) indicate that the potential $q_{i}(x)$ can be
recovered from the first coefficients of the series (\ref{phiNSBF}) or
(\ref{S}). Indeed, we have
\begin{equation}
q_{i}(x)=\frac{g_{i,0}^{\prime\prime}(x)}{g_{i,0}(x)+1} \label{qi from g0}%
\end{equation}
and
\begin{equation}
q_{i}(x)=\frac{\left(  xs_{i,0}(x)\right)  ^{\prime\prime}}{xs_{i,0}(x)+3x}.
\label{qi from s0}%
\end{equation}

\section{Leaf peeling\label{Sect LeafPeeling}}

The overall strategy is based on the leaf peeling method (developed in
\cite{AvdoninKurasov2008}) combined with the use of the NSBF representations
from Section \ref{Sect FundSys}.

\begin{definition}
A star shaped subgraph of a tree graph is called a \textbf{sheaf} if all but
one of its edges are leaf edges. The only internal vertex of the sheaf is
called the \textbf{abscission vertex}. The only edge which is not a leaf edge
is called the \textbf{stem edge }of the sheaf.
\end{definition}

\begin{figure}
[ptb]
\begin{center}
\includegraphics[
height=2.5482in,
width=4.1424in
]%
{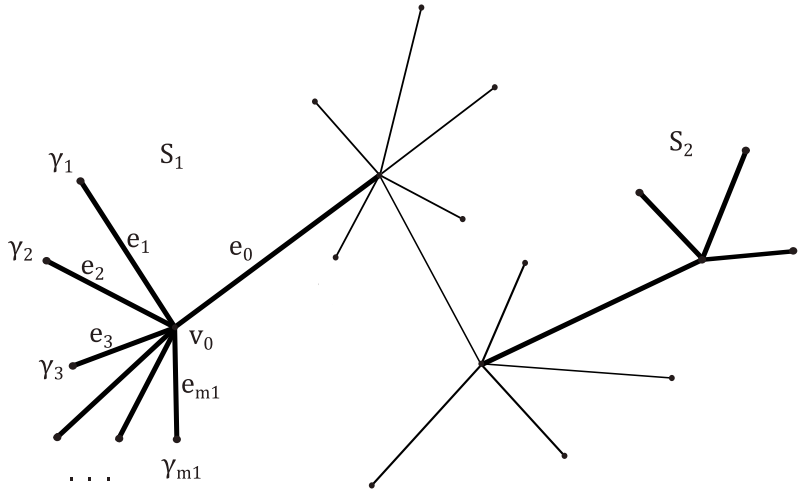}%
\caption{An example of a tree graph with two sheaves $S_{1}$ and $S_{2}$
highlighted in bold lines. Here $v_{0}$ is the abscission vertex and $e_{0}$
is the stem edge of $S_{1}$.}%
\label{FigGraph}%
\end{center}
\end{figure}

In Fig. \ref{FigGraph} an example of a tree graph with exactly two sheaves
$S_{1}$ and $S_{2}$ is shown. Here $v_{0}$ is the abscission vertex and
$e_{0}$ is the stem edge of $S_{1}$.%

A tree graph, which is not a star shaped graph, contains at least one sheaf
\cite{AChoque}, \cite{AvdoninZhao2020}. The method developed in \cite{AKK2023}
allows one to recover the potential on all leaf edges of a sheaf from several
entries of the Weyl matrix that correspond to the leaf edges. After that, the
leaf peeling method serves for calculating the Weyl matrix for the new
(smaller) tree graph obtained by cutting off the leaf edges of the sheaf, so
that the abscission vertex becomes a new leaf.

Repeating this leaf peeling procedure one comes eventually to an inverse
problem of recovering the potential of a star shaped graph from its Weyl
matrix. For this problem the method based on the NSBF representations was
developed in \cite{AKK2023}. First of all, let us explain the leaf peeling
method in more detail.

Consider a tree graph $\Omega$, and assume $S_{1}$ to be its sheaf with an
abscission vertex $v_{0}$ and a stem edge $e_{0}$. Without loss of generality
we assume that $\gamma_{1},\ldots,\gamma_{m_{1}}$ are the leaves belonging to
$S_{1}$. By $\widetilde{\Omega}$ we denote the subgraph $\widetilde{\Omega
}=\Omega\setminus\left\{  e_{1},\ldots,e_{m_{1}}\right\}  $, which is obtained
by removing the leaf edges of $S_{1}$ from $\Omega$, so that the vertex
$v_{0}$ is a leaf of $\widetilde{\Omega}$, and $e_{0}$ is the corresponding
leaf edge. We write $\gamma_{0}=v_{0}$. By $\widetilde{\mathbf{M}}$ we denote
the Weyl matrix of $\widetilde{\Omega}$.

Assume that the potential $q(x)$ is already known on each\ leaf edge
$e_{1},\ldots,e_{m_{1}}$ of the sheaf $S_{1}$, so that $q_{1}(x),\ldots
,q_{m_{1}}(x)$ are known as well as
\begin{equation}
\varphi_{i}(\rho,L_{i}),\,S_{i}(\rho,L_{i}),\,\varphi_{i}^{\prime}(\rho
,L_{i}),\,S_{i}^{\prime}(\rho,L_{i})\quad\text{for\thinspace}1\leq i\leq
m_{1}. \label{phiS}%
\end{equation}
The problem which is solved by the leaf peeling method consists in computing
$\widetilde{\mathbf{M}}$ from $\mathbf{M}$ and (\ref{phiS}).

We begin with the element $\widetilde{\mathbf{M}}_{00}(\rho^{2})$, which is
the derivative at $\gamma_{0}$ of the Weyl solution $\widetilde{w}_{0}%
(\rho,x)$ on $\widetilde{\Omega}$, which equals one at $\gamma_{0}$ and zero
at all other leaves of $\widetilde{\Omega}$, that is, at $\gamma_{j}$ for
$j=m_{1}+1,\ldots,m$. Note that
\[
\widetilde{w}_{0}(\rho,x)=\frac{w_{1}(\rho,x)}{w_{11}(\rho,\gamma_{0})}%
\quad\text{on }\widetilde{\Omega}.
\]
Hereafter, the notation $w_{ij}(\rho,x)$ means that we consider $j$-th
component of a solution $w_{i}(\rho,x)$, that is, the solution $w_{i}(\rho,x)$
on the edge $e_{j}$. Thus,
\[
\widetilde{\mathbf{M}}_{00}(\rho^{2})=\frac{\partial w_{10}(\rho,\gamma_{0}%
)}{w_{11}(\rho,\gamma_{0})}=-\sum_{j=1}^{m_{1}}\frac{\partial w_{1j}%
(\rho,\gamma_{0})}{w_{11}(\rho,\gamma_{0})},
\]
where $w_{10}(\rho,x)$ is a component of $w_{1}(\rho,x)$ on $e_{0}$. This can
be written in the form%
\[
\widetilde{\mathbf{M}}_{00}(\rho^{2})=
\]%
\[
-\frac{1}{\varphi_{1}(\rho,L_{1})+\mathbf{M}_{11}(\rho^{2})S_{1}(\rho,L_{1}%
)}\left(  \varphi_{1}^{\prime}(\rho,L_{1})+\mathbf{M}_{11}(\rho^{2}%
)S_{1}^{\prime}(\rho,L_{1})+\sum_{j=2}^{m_{1}}\mathbf{M}_{1j}(\rho^{2}%
)S_{j}^{\prime}(\rho,L_{j})\right)  .
\]

Next, it is immediate to see that
\[
\widetilde{\mathbf{M}}_{0i}(\rho^{2})=\partial\widetilde{w}_{0}(\rho
,\gamma_{i})=\frac{\partial w_{1i}(\rho,\gamma_{i})}{w_{11}(\rho,\gamma_{0}%
)}=\frac{\mathbf{M}_{1i}(\rho^{2})}{w_{11}(\rho,\gamma_{0})}=\frac
{\mathbf{M}_{1i}(\rho^{2})}{\varphi_{1}(\rho,L_{1})+\mathbf{M}_{11}(\rho
^{2})S_{1}(\rho,L_{1})}%
\]
for $i=m_{1}+1,\ldots,m$. Thus, we obtained the first row of the Weyl matrix
$\widetilde{\mathbf{M}}(\rho^{2})$.

Now choose $i\in\left\{  m_{1}+1,\ldots,m\right\}  $ and consider the
following solution on $\Omega$:%
\[
u(\rho,x):=w_{i}(\rho,x)-\alpha w_{1}(\rho,x)
\]
where $\alpha$ is a constant to be defined. Note that%
\[
u(\rho,\gamma_{i})=1
\]
and
\[
u(\rho,\gamma_{j})=0\quad\text{for\thinspace}j=m_{1}+1,\ldots,m\text{ and
}j\neq i.
\]

Moreover, $\alpha$ can be chosen such that additionally $u(\rho,\gamma_{0}%
)=0$. Indeed, take
\[
\alpha=\frac{w_{i}(\rho,\gamma_{0})}{w_{1}(\rho,\gamma_{0})}.
\]
Then
\[
u(\rho,\gamma_{0})=w_{i}(\rho,\gamma_{0})-\frac{w_{i}(\rho,\gamma_{0})}%
{w_{1}(\rho,\gamma_{0})}w_{1}(\rho,\gamma_{0})=0
\]
and hence with this choice of $\alpha$ we have that
\begin{equation}
\widetilde{w}_{i}(\rho,x)=u(\rho,x)=w_{i}(\rho,x)-\frac{w_{i}(\rho,\gamma
_{0})}{w_{1}(\rho,\gamma_{0})}w_{1}(\rho,x)\quad\text{for\thinspace}%
i=m_{1}+1,\ldots,m \label{Weyl solution new}%
\end{equation}
is the Weyl solution on $\widetilde{\Omega}$. Thus,
\begin{align}
\widetilde{\mathbf{M}}_{i0}(\rho^{2})  &  =\widetilde{w}_{i0}^{\prime}%
(\rho,\gamma_{0})=w_{i0}^{\prime}(\rho,\gamma_{0})-\frac{w_{i1}(\rho
,\gamma_{0})}{w_{11}(\rho,\gamma_{0})}w_{10}^{\prime}(\rho,\gamma
_{0})\nonumber\\
&  =w_{i0}^{\prime}(\rho,\gamma_{0})-w_{i1}(\rho,\gamma_{0})\widetilde
{\mathbf{M}}_{00}(\rho^{2}). \label{Mtili0}%
\end{align}

We have
\[
w_{i1}(\rho,\gamma_{0})=\mathbf{M}_{i1}(\rho^{2})S_{1}(\rho,L_{1})
\]
and
\[
w_{i0}^{\prime}(\rho,\gamma_{0})=-\sum_{j=1}^{m_{1}}\mathbf{M}_{ij}(\rho
^{2})S_{j}^{\prime}(\rho,L_{j}).
\]
Thus, all the terms in (\ref{Mtili0}) are known, so that the first column
(corresponding to $\gamma_{0}$) of the Weyl matrix $\widetilde{\mathbf{M}%
}(\rho^{2})$ can be computed.

Finally, for $j=m_{1}+1,\ldots,m$ from (\ref{Weyl solution new}) we have
\begin{align*}
\widetilde{\mathbf{M}}_{ij}(\rho^{2})  &  =\partial\widetilde{w}_{i}%
(\rho,\gamma_{j})=\mathbf{M}_{ij}(\rho^{2})-\frac{w_{i}(\rho,\gamma_{0}%
)}{w_{1}(\rho,\gamma_{0})}\mathbf{M}_{1j}(\rho^{2})\\
&  =\mathbf{M}_{ij}(\rho^{2})-w_{i}(\rho,\gamma_{0})\widetilde{\mathbf{M}%
}_{0j}(\rho^{2})\\
&  =\mathbf{M}_{ij}(\rho^{2})-\mathbf{M}_{i1}(\rho^{2})S_{1}(\rho
,L_{1})\widetilde{\mathbf{M}}_{0j}(\rho^{2}).
\end{align*}
This finishes the calculation of the Weyl matrix $\widetilde{\mathbf{M}}%
(\rho^{2})$.

\section{Solution of local problem\label{Sect Local Problem}}

By the \textbf{local problem }we understand the recovery of the potential
$q(x)$ on the leaf edges $e_{1},\ldots,e_{m_{1}}$ of a sheaf, and moreover of
all the functions in (\ref{phiS}).

As we see from the previous section, the leaf peeling method gives us the
possibility to reduce the solution of the inverse problem on the tree graph to
a successive solution of local problems on sheaves. After having solved each
local problem, the sheaf can be removed, and the Weyl matrix for the remaining
tree graph is calculated using the leaf peeling method. In the last step,
after having removed all possible sheaves, one obtains an inverse problem on a
star shaped graph, which in fact is solved by the same method as we explain in
the present section for local problems.

Let $S$ be a sheaf, whose leaf edges are $e_{1},\ldots,e_{m_{1}}$. We assume
that $\mathbf{M}_{ij}(\rho^{2})$ are known for $1\leq i,j\leq m_{1}$ and on a
set of values $\rho=\rho_{k}$. From this information we look to recover the
potential $q_{i}(x)$, $i=1,\ldots,m_{1}$. The proposed method for solving this
problem includes several steps.

First, we compute a number of the constants $\left\{  g_{i,n}(L_{i})\right\}
_{n=0}^{N}$ and $\left\{  s_{i,n}(L_{i})\right\}  _{n=0}^{N}$ for every
$i=1,\ldots,m_{1}$, that is, the values of the coefficients from
(\ref{phiNSBF}) and (\ref{S}) at the endpoint of the edge $e_{i}$. This first
step allows us to split the problem reducing it to separate inverse problems
on the edges. Indeed, the knowledge of the coefficients $\left\{
g_{i,n}(L_{i})\right\}  _{n=0}^{N}$ and $\left\{  s_{i,n}(L_{i})\right\}
_{n=0}^{N}$ allows us in the second step to compute the Dirichlet-Dirichlet
and Neumann-Dirichlet spectra for the potential $q_{i}(x)$, $x\in\left[
0,L_{i}\right]  $, thus obtaining a two-spectra inverse problem for $q_{i}(x)$.

In the third step the two-spectra inverse problem is solved with the aid of
the representation (\ref{phiNSBF}) and an analogous series representation for
the solution $T_{i}(\rho,x)$ of (\ref{Schri}) satisfying the initial
conditions at the endpoint $L_{i}$:
\begin{equation}
T_{i}(\rho,L_{i})=0,\quad T_{i}^{\prime}(\rho,L_{i})=1. \label{psi init}%
\end{equation}
The two-spectra inverse problem is reduced to a system of linear algebraic
equations, from which we obtain the coefficient $g_{i,0}(x)$. Finally, the
potential $q_{i}(x)$ is calculated from (\ref{qi from g0}).

\subsection{Calculation of coefficients $\left\{  g_{i,n}(L_{i})\right\}
_{n=0}^{N}$ and $\left\{  s_{i,n}(L_{i})\right\}  _{n=0}^{N}$
\label{Subsect Coefs}}

The knowledge of $\mathbf{M}_{ij}(\rho_{k}^{2})$, $1\leq i,j\leq m_{1}$ means
that for each $i=1,\ldots,m_{1}$ we have the equalities
\[
\varphi_{i}(\rho_{k},L_{i})+\mathbf{M}_{ii}(\rho_{k}^{2})S_{i}(\rho_{k}%
,L_{i})=\mathbf{M}_{ij}(\rho_{k}^{2})S_{j}(\rho_{k},L_{j}),
\]%
\[
\mathbf{M}_{ij}(\rho_{k}^{2})S_{j}(\rho_{k},L_{j})=\mathbf{M}_{il}(\rho
_{k}^{2})S_{l}(\rho_{k},L_{l}),
\]
valid for all $j,l\neq i$, $1\leq j,l\leq m_{1}$. Thus, for every $\rho_{k}$,
using the representations (\ref{phiNSBF}) and (\ref{S}) we have the equations
\begin{align}
&  \rho_{k}\sum_{n=0}^{\infty}(-1)^{n}g_{i,n}(L_{i})\mathbf{j}_{2n}(\rho
_{k}L_{i})+\mathbf{M}_{ii}(\rho_{k}^{2})\sum_{n=0}^{\infty}(-1)^{n}%
s_{i,n}(L_{i})\mathbf{j}_{2n+1}(\rho_{k}L_{i})\nonumber\\
&  -\mathbf{M}_{ii+1}(\rho_{k}^{2})\sum_{n=0}^{\infty}(-1)^{n}s_{i+1,n}%
(L_{i+1})\mathbf{j}_{2n+1}(\rho_{k}L_{i+1})\nonumber\\
&  =\mathbf{M}_{ii+1}(\rho_{k}^{2})\sin(\rho_{k}L_{i+1})-\rho_{k}\cos(\rho
_{k}L_{i})-\mathbf{M}_{ii}(\rho_{k}^{2})\sin(\rho_{k}L_{i}),\quad\text{for
}i=1,\ldots,m_{1}, \label{type1}%
\end{align}
where for $i=m_{1}$ we replace $i+1$ by $1$ (the cyclic rule), as well as the
equations%
\begin{align}
&  \mathbf{M}_{ij}(\rho_{k}^{2})\sum_{n=0}^{\infty}(-1)^{n}s_{j,n}%
(L_{j})\mathbf{j}_{2n+1}(\rho_{k}L_{j})-\mathbf{M}_{ij+1}(\rho_{k}^{2}%
)\sum_{n=0}^{\infty}(-1)^{n}s_{j+1,n}(L_{j+1})\mathbf{j}_{2n+1}(\rho
_{k}L_{j+1})\nonumber\\
&  =\mathbf{M}_{ij+1}(\rho_{k}^{2})\sin(\rho_{k}L_{j+1})-\mathbf{M}_{ij}%
(\rho_{k}^{2})\sin(\rho_{k}L_{j}), \label{type2}%
\end{align}
where $j=1,\ldots,m_{1}$, $j\neq i$, $j+1\neq i$, and again, for $j=m_{1}$ we
replace $j+1$ by $1$.

Suppose that $\mathbf{M}_{ij}(\rho_{k}^{2})$, $1\leq i,j\leq m_{1}$ are given
at $K$ points $\rho_{k}^{2}$. To compute the finite sets of the coefficients
$\left\{  g_{i,n}(L_{i})\right\}  _{n=0}^{N}$ and $\left\{  s_{i,n}%
(L_{i})\right\}  _{n=0}^{N}$ for all $i=1,\ldots,m_{1}$ we chose the following
strategy. Fix $i$ and consider corresponding equations (\ref{type1}) and
(\ref{type2}). For each $i$ we have one equation of the form (\ref{type1}) and
$m_{1}-2$ equations of the form (\ref{type2}). Considering these $m_{1}-1$
equations for every $\rho_{k}^{2}$, we obtain $K(m_{1}-1)$ equations for
$(m_{1}+1)(N+1)$ unknowns. Here the unknowns are $m_{1}$ sets of the
coefficients $\left\{  s_{j,n}(L_{i})\right\}  _{n=0}^{N}$ (for $j=1,\ldots
,m_{1}$) and one set of the coefficients $\left\{  g_{i,n}(L_{i})\right\}
_{n=0}^{N}$. Thus, we need at least $K=\left\lceil \frac{(m_{1}+1)(N+1)}%
{m_{1}-1}\right\rceil $ points $\rho_{k}^{2}$ at which the elements of the
Weyl matrix $\mathbf{M}_{ij}(\rho_{k}^{2})$, $1\leq i,j\leq m_{1}$ are known.
Here $\left\lceil x\right\rceil $ denotes the least integer greater than or
equal to $x$. Thus obtained system of equations allows us to compute $\left\{
g_{i,n}(L_{i})\right\}  _{n=0}^{N}$ and $\left\{  s_{i,n}(L_{i})\right\}
_{n=0}^{N}$ for the fixed $i$. In total, we consider $m_{1}$ linear algebraic
systems of this kind to compute all sets of the coefficients $\left\{
g_{i,n}(L_{i})\right\}  _{n=0}^{N}$ and $\left\{  s_{i,n}(L_{i})\right\}
_{n=0}^{N}$ for all $i=1,\ldots,m_{1}$.

An important observation however consists in the fact that in order to find
the coefficients $\left\{  g_{i,n}(L_{i})\right\}  _{n=0}^{N}$ and $\left\{
s_{i,n}(L_{i})\right\}  _{n=0}^{N}$, $i=1,\ldots,m_{1}$, there is no need to
consider the equations of the form (\ref{type2}) or at least all such
equations. In fact, it is enough to consider equation (\ref{type1}) alone,
which means that to recover the potential on the edge $e_{i}$, we need to know
the main diagonal entry \textbf{$M$}$_{ii}(\rho_{k}^{2})$ as well as some
$\mathbf{M}_{ij}(\rho_{k}^{2})$ for one $j\neq i$, which can be $\mathbf{M}%
_{ii+1}(\rho_{k}^{2})$ (for $i=m_{1}$, $i+1$ is replaced by $1$). In this
case, for each edge $e_{i}$, from equations of the form (\ref{type1}), which
we have for each $\rho_{k}$, we compute $\left\{  g_{i,n}(L_{i})\right\}
_{n=0}^{N}$, $\left\{  s_{i,n}(L_{i})\right\}  _{n=0}^{N}$ and $\left\{
s_{i+1,n}(L_{i+1})\right\}  _{n=0}^{N}$, that is, $3(N+1)$ unknowns. As it was
shown in \cite{AKK2023}, eventually the accuracy of the recovered potential is
comparable when one uses all equations of the form (\ref{type2}), part of them
or even none.

\subsection{Reduction to the two-spectra inverse problem on the edge
\label{Subsect reduction to two spectra}}

The first step, as described in the previous subsection, reduces the problem
on the sheaf to $m_{1}$ separate problems on the edges. Thus, consider an edge
$e_{i}$ for which at this stage we have computed the coefficients $\left\{
g_{i,n}(L_{i})\right\}  _{n=0}^{N}$ and $\left\{  s_{i,n}(L_{i})\right\}
_{n=0}^{N}$. Now we use them to compute the Dirichlet-Dirichlet and
Neumann-Dirichlet spectra for the potential $q_{i}(x)$, $x\in\left[
0,L_{i}\right]  $. This is done with the aid of the approximate solutions
evaluated at the end point%
\[
\varphi_{i,N}(\rho,L_{i})=\cos\left(  \rho L_{i}\right)  +\sum_{n=0}%
^{N}(-1)^{n}g_{i,n}(L_{i})\mathbf{j}_{2n}(\rho L_{i})
\]
and
\[
S_{i,N}(\rho,L_{i})=\frac{\sin\left(  \rho L_{i}\right)  }{\rho}+\frac{1}%
{\rho}\sum_{n=0}^{N}(-1)^{n}s_{i,n}(L_{i})\mathbf{j}_{2n+1}(\rho L_{i}).
\]
Indeed, since $S_{i}(\rho,x)$ (as well as $S_{i,N}(\rho,x)$) satisfies the
Dirichlet condition at the origin, zeros of $S_{i}(\rho,L_{i})$ are precisely
square roots of the Dirichlet-Dirichlet eigenvalues. That is, $S_{i}%
(\rho,L_{i})$ is the characteristic function of the Sturm-Liouville problem
\begin{equation}
-y^{\prime\prime}+q_{i}(x)y=\lambda y,\quad x\in(0,L_{i}), \label{SLi}%
\end{equation}%
\begin{equation}
y(0)=y(L_{i})=0, \label{DD cond}%
\end{equation}
and its zeros coincide with the numbers $\left\{  \mu_{i,k}\right\}
_{k=1}^{\infty}$, such that $\mu_{i,k}^{2}$ are the eigenvalues of the problem
(\ref{SLi}), (\ref{DD cond}). In turn, zeros of the function $S_{i,N}%
(\rho,L_{i})$ approximate zeros of $S_{i}(\rho,L_{i})$ (see Theorem
\ref{Th closeness of zeros} above). Thus, the singular numbers $\left\{
\mu_{i,k}\right\}  _{k=1}^{\infty}$ are approximated by zeros of the function
$S_{i,N}(\rho,L_{i})$.

The same reasoning is valid for the function $\varphi_{i,N}(\rho,L_{i})$,
whose zeros approximate the singular numbers $\left\{  \nu_{i,k}\right\}
_{k=1}^{\infty}$, which are the square roots of the Neumann-Dirichlet
eigenvalues,\ i.e., the eigenvalues of the Sturm-Liouville problem for
(\ref{SLi}) subject to the boundary conditions
\begin{equation}
y^{\prime}(0)=y(L_{i})=0. \label{ND cond}%
\end{equation}
\qquad

Thus, on every edge $e_{i}$ we obtain the classical inverse problem of
recovering the potential $q_{i}(x)$ from two spectra, which is considered in
the next step.

\subsection{Solution of two-spectra inverse
problem\label{Subsect Solution two-spectra}}

For solving the obtained inverse problem we use the method developed in
\cite{AKK2023}. For the sake of completeness we briefly describe it here. At
this stage we dispose of two finite sequences of singular numbers $\left\{
\mu_{i,k}\right\}  _{k=1}^{K_{D}}$ and $\left\{  \nu_{i,k}\right\}
_{k=1}^{K_{N}}$ which are square roots of the eigenvalues of problems
(\ref{SLi}), (\ref{DD cond}) and (\ref{SLi}), (\ref{ND cond}), respectively,
as well as of two sequences of numbers $\left\{  s_{i,n}(L_{i})\right\}
_{n=0}^{N}$ and $\left\{  g_{i,n}(L_{i})\right\}  _{n=0}^{N}$, which are the
values of the coefficients from (\ref{S}) and (\ref{phiNSBF}) at the endpoint.

Let us consider the solution $T_{i}(\rho,x)$ of equation (\ref{SLi})
satisfying the initial conditions at $L_{i}$:
\[
T_{i}(\rho,L_{i})=0,\quad T_{i}^{\prime}(\rho,L_{i})=1.
\]
Analogously to the solution (\ref{S}), the solution $T_{i}(\rho,x)$ admits the
series representation%
\begin{equation}
T_{i}(\rho,x)=\frac{\sin\left(  \rho\left(  x-L_{i}\right)  \right)  }{\rho
}+\frac{1}{\rho}\sum_{n=0}^{\infty}(-1)^{n}t_{i,n}(x)\mathbf{j}_{2n+1}%
(\rho\left(  x-L_{i}\right)  ), \label{Ti}%
\end{equation}
where $t_{i,n}\left(  x\right)  $ are corresponding coefficients, analogous to
$s_{i,n}\left(  x\right)  $ from (\ref{S}).

Note that for $\rho=\nu_{i,k}$ the solutions $\varphi_{i}(\nu_{i,k},x)$ and
$T_{i}(\nu_{i,k},x)$ are linearly dependent because both are eigenfunctions of
problem (\ref{SLi}), (\ref{ND cond}). Hence there exist such real constants
$\beta_{i,k}\neq0$, that
\begin{equation}
\varphi_{i}(\nu_{i,k},x)=\beta_{i,k}T_{i}(\nu_{i,k},x). \label{phi=T}%
\end{equation}
Moreover, these multiplier constants can be easily calculated by recalling
that $\varphi_{i}(\nu_{i,k},0)=1$. Thus,%
\begin{align}
\frac{1}{\beta_{i,k}}  &  =T_{i}(\nu_{i,k},0)\approx T_{i,N}(\nu
_{i,k},0)\nonumber\\
&  =-\frac{\sin\left(  \nu_{i,k}L_{i}\right)  }{\nu_{i,k}}-\frac{1}{\nu_{i,k}%
}\sum_{n=0}^{N}(-1)^{n}t_{i,n}(0)\mathbf{j}_{2n+1}(\nu_{i,k}L_{i}),
\label{betaik}%
\end{align}
where we took into account that the spherical Bessel functions of odd order
are odd functions. The coefficients $\left\{  t_{i,n}(0)\right\}  _{n=0}^{N}$
are computed with the aid of the singular numbers $\left\{  \mu_{i,k}\right\}
_{k=1}^{K_{D}}$ as follows. Since the functions $T_{i}(\mu_{i,k},x)$,
$k=1,2,\ldots$ are eigenfunctions of the problem (\ref{SLi}), (\ref{DD cond}),
we have that $T_{i}(\mu_{i,k},0)=0$ and hence%
\[
\sum_{n=0}^{\infty}(-1)^{n}t_{i,n}(0)\mathbf{j}_{2n+1}(\mu_{i,k}L_{i}%
)=-\sin\left(  \mu_{i,k}L_{i}\right)  ,\quad k=1,2,\cdots.
\]
This leads to a system of linear algebraic equations for computing the
coefficients $\left\{  t_{i,n}(0)\right\}  _{n=0}^{N}$, which has the form%
\[
\sum_{n=0}^{N}(-1)^{n}t_{i,n}(0)\mathbf{j}_{2n+1}(\mu_{i,k}L_{i})=-\sin\left(
\mu_{i,k}L_{i}\right)  ,\quad k=1,\cdots,K_{D}.
\]
Now, having computed $\left\{  t_{i,n}(0)\right\}  _{n=0}^{N}$, we compute the
multiplier constants $\left\{  \beta_{i,k}\right\}  _{k=1}^{K_{N}}$ from
(\ref{betaik}).

Next, we use equation (\ref{phi=T}) for constructing a system of linear
algebraic equations for the coefficients $g_{i,n}(x)$ and $t_{i,n}\left(
x\right)  $. Indeed, equation (\ref{phi=T}) can be written in the form%
\begin{align*}
&  \sum_{n=0}^{\infty}(-1)^{n}g_{i,n}(x)\mathbf{j}_{2n}(\nu_{i,k}%
x)-\frac{\beta_{i,k}}{\nu_{i,k}}\sum_{n=0}^{\infty}(-1)^{n}t_{i,n}%
(x)\mathbf{j}_{2n+1}(\nu_{i,k}\left(  x-L_{i}\right)  )\\
&  =\frac{\beta_{i,k}}{\nu_{i,k}}\sin\left(  \nu_{i,k}\left(  x-L_{i}\right)
\right)  -\cos\left(  \nu_{i,k}x\right)  .
\end{align*}

We have as many of such equations as many Neumann-Dirichlet singular numbers
$\nu_{i,k}$ are computed. For computational purposes we choose some natural
number $N_{c}$ - the number of the coefficients $g_{i,n}(x)$ and
$t_{i,n}\left(  x\right)  $ to be computed. More precisely, we choose a
sufficiently dense set of points $x_{m}\in(0,L_{i})$ and at every $x_{m}$
consider the equations%
\begin{align*}
&  \sum_{n=0}^{N_{c}}(-1)^{n}g_{i,n}(x_{m})\mathbf{j}_{2n}(\nu_{i,k}%
x_{m})-\frac{\beta_{i,k}}{\nu_{i,k}}\sum_{n=0}^{N_{c}}(-1)^{n}t_{i,n}%
(x_{m})\mathbf{j}_{2n+1}(\nu_{i,k}\left(  x_{m}-L_{i}\right)  )\\
&  =\frac{\beta_{i,k}}{\nu_{i,k}}\sin\left(  \nu_{i,k}\left(  x_{m}%
-L_{i}\right)  \right)  -\cos\left(  \nu_{i,k}x_{m}\right)  ,\quad
k=1,\ldots,K_{N}.
\end{align*}
Solving this system of equations we find $g_{i,0}(x_{m})$ and consequently
$g_{i,0}(x)$ at a sufficiently dense set of points of the interval $(0,L_{i})$
as well as $t_{i,0}(x)$ (which is used below). Finally, with the aid of
(\ref{qi from g0}) we compute $q_{i}(x)$.

\subsection{Computation of functions (\ref{phiS}%
)\label{Subsect computation of phi S}}

Having computed $q_{i}(x)$ on each leaf edge $e_{i}$ of the sheaf, in
principle, gives us the possibility to compute the functions (\ref{phiS})
(needed for applying the leaf peeling procedure) by solving corresponding
Cauchy problems on $e_{i}$. However, this is clearly not the most attractive
option because it implies solving the Cauchy problems for $K$ different values
of $\rho_{k}$ and for numerically recovered $q_{i}(x)$. Instead, we use again
the NSBF representations and the fact that in the preceding step we computed
the functions $g_{i,0}(x)$ and $t_{i,0}(x)$. In \cite{KNT} (see also
\cite[Sect. 9.4]{KrBook2020}) a recurrent integration procedure was developed
for calculating the coefficients of the NSBF representations (\ref{phiNSBF}%
)-(\ref{Sprime}). It requires the knowledge of a nonvanishing on $\left[
0,L_{i}\right]  $ solution $f(x)$ of the equation
\[
f^{\prime\prime}(x)-q_{i}(x)f(x)=0
\]
satisfying the initial condition
\[
f(0)=1,
\]
as well as of its derivative $f^{\prime}(x)$. Denote
\[
h:=f^{\prime}(0),
\]
which is a complex number. Then all the NSBF coefficients in (\ref{phiNSBF}%
)-(\ref{Sprime}) can be calculated from $f(x)$ and $f^{\prime}(x)$ with the
aid of a recurrent integration procedure. More precisely, we will calculate
the solutions $\varphi_{i,h}(\rho,x)$, $S_{i}(\rho,x)$ and their first
derivatives, where $\varphi_{i,h}(\rho,x)$ is a solution of (\ref{Schri})
satisfying the initial conditions
\[
\varphi_{i,h}(\rho,0)=1,\quad\varphi_{i,h}^{\prime}(\rho,0)=h,
\]
and $S_{i}(\rho,x)$ is the solution defined above. Obviously,
\begin{equation}
\varphi_{i}(\rho,x)=\varphi_{i,h}(\rho,x)-hS_{i}(\rho,x). \label{phi_i=}%
\end{equation}

First, let us obtain such a nonvanishing solution $f(x)$ and then briefly
remind the recurrent integration procedure.

From (\ref{beta0}) we have that
\[
\varphi_{i}(0,x)=g_{i,0}(x)+1.
\]
Also, from (\ref{Ti}) we have
\[
T_{i}(0,x)=\left(  x-L_{i}\right)  \left(  \frac{t_{i,0}(x)}{3}+1\right)  ,
\]
where we took into account that $\mathbf{j}_{1}(z)\sim\frac{z}{3}$,
$z\rightarrow0$. Assume first that zero is not a Neumann-Dirichlet eigenvalue,
i.e., $\varphi_{i}(0,L_{i})\neq0$. Then $\varphi_{i}(0,x)$ and $T_{i}(0,x)$
are linearly independent, and hence the complex valued solution
\[
f(x):=\varphi_{i}(0,x)+\mathbf{i}T_{i}(0,x)
\]
has no zero in $\left[  0,L_{i}\right]  $ ($\mathbf{i}$ stands for the
imaginary unit). We have $h=f^{\prime}(0)=\mathbf{i}T_{i}^{\prime
}(0,0)=\mathbf{i}\left(  1+\frac{1-L_{i}}{3}t_{i,0}^{\prime}(0)\right)  $. If,
on the contrary, zero is a Neumann-Dirichlet eigenvalue, we can construct the
linearly independent solution using the Abel formula $\psi_{i}(x):=\varphi
_{i}(0,x)\int_{0}^{x}\frac{dt}{\varphi_{i}^{2}(0,x)}$. Then $f(x):=\varphi
_{i}(0,x)+\mathbf{i}\psi_{i}(x)$ and $h=-\frac{\mathbf{i}}{\varphi_{i}%
^{\prime}(0,L_{i})}$. In any case, we have a nonvanishing solution $f(x)$ and
can apply the recurrent integration procedure. With its aid we compute two
sequences of functions $\left\{  \beta_{n}(x)\right\}  _{n=-1}^{\infty}$ and
$\left\{  \xi_{n}(x)\right\}  _{n=-1}^{\infty}$. The first elements are
defined by the equalities
\[
\beta_{-1}(x)=\frac{1}{2},\quad\beta_{0}(x)=\frac{(f(x)-1)}{2},\quad\xi
_{-1}(x)=\frac{1}{4}\int_{0}^{x}q(s)\,ds,\quad\xi_{0}(x)=\frac{f^{\prime
}(x)-h}{2}-\frac{1}{4}\int_{0}^{x}q(s)\,ds,
\]
while all other elements are constructed as follows%
\[
\beta_{n}(x)=\frac{2n+1}{2n-3}\left(  \beta_{n-2}(x)+\frac{c_{n}}{x^{n}%
}f(x)\theta_{n}(x)\right)  ,
\]%
\[
\xi_{n}(x)=\frac{2n+1}{2n-3}\left(  \xi_{n-2}(x)+\frac{c_{n}}{x^{n}}\left(
f^{\prime}(x)\theta_{n}(x)+\frac{\eta_{n}(x)}{f(x)}\right)  -\frac{c_{n}%
-2n+1}{x}\beta_{n-2}(x)\right)  ,
\]
where
\[
\eta_{n}(x)=\int_{0}^{x}\bigl(tf^{\prime}(t)+(n-1)f(t)\bigr)\beta
_{n-2}(t)t^{n-2}\,dt,\quad\theta_{n}(x)=\int_{0}^{x}\frac{1}{f^{2}%
(t)}\bigl(\eta_{n}(t)-f(t)\beta_{n-2}(t)t^{n-1}\bigr)dt
\]
for $n=1,2,\ldots$, and $c_{n}=1$ if $n=1$ and $c_{n}=2(2n-1)$ otherwise.

Finally,
\[
\varphi_{i,h}(\rho,x)=\cos\left(  \rho x\right)  +2\sum_{n=0}^{\infty}%
(-1)^{n}\beta_{2n}(x)\mathbf{j}_{2n}(\rho x),
\]%
\[
S_{i}(\rho,x)=\frac{\sin\left(  \rho x\right)  }{\rho}+\frac{2}{\rho}%
\sum_{n=0}^{\infty}(-1)^{n}\beta_{2n+1}(x)\mathbf{j}_{2n+1}(\rho x),
\]%
\[
\varphi_{i,h}^{\prime}(\rho,x)=-\rho\sin\left(  \rho x\right)  +\cos\left(
\rho x\right)  \left(  \frac{1}{2}\int_{0}^{x}q_{i}(t)\,dt+h\right)
+2\sum_{n=0}^{\infty}(-1)^{n}\xi_{2n}(x)\mathbf{j}_{2n}(\rho x)
\]
and
\[
S_{i}^{\prime}(\rho,x)=\cos\left(  \rho x\right)  +\frac{\sin\left(  \rho
x\right)  }{2\rho}\int_{0}^{x}q_{i}(t)\,dt+\frac{2}{\rho}\sum_{n=0}^{\infty
}(-1)^{n}\xi_{2n+1}(x)\mathbf{j}_{2n+1}(\rho x).
\]
These formulas together with (\ref{phi_i=}) lead to the computation of the
functions (\ref{phiS}).

\section{Summary of the method\label{Sect SummaryMethod}}

Given a tree graph $\Omega$ and the Weyl matrix $\mathbf{M}(\rho^{2})$ for a
number of values $\rho_{k}^{2}$, $k=1,\ldots,K$. In the first step we identify
a sheaf $S_{1}$ and apply to it the procedure for solving the local inverse
problem, explained in Section \ref{Sect Local Problem}. According to it, for
each leaf edge $e_{i}$ of the sheaf we compute the sets of the NSBF
coefficients $\left\{  g_{i,n}(L_{i})\right\}  _{n=0}^{N}$ and $\left\{
s_{i,n}(L_{i})\right\}  _{n=0}^{N}$ (see subsection \ref{Subsect Coefs}).
Next, we compute the Dirichlet-Dirichlet and Neumann-Dirichlet spectra on each
leaf edge (see subsection \ref{Subsect reduction to two spectra}), thus
obtaining an inverse two-spectra problem on $e_{i}$. This problem is solved by
the method described in subsection \ref{Subsect Solution two-spectra}. We
obtain $q_{i}(x)$ and additionally the functions (\ref{phiS}), as explained in
subsection \ref{Subsect computation of phi S}. This gives us a complete
solution of the local problem on the sheaf $S_{1}$. In the next step the leaf
edges of $S_{1}$ are removed, and the Weyl matrix $\widetilde{\mathbf{M}}%
(\rho^{2})$ for the new smaller tree graph $\widetilde{\Omega}$ is calculated
following the leaf peeling method from Section \ref{Sect LeafPeeling}.

This sequence of steps is repeated until arriving at a final star shaped
graph, for which the corresponding Weyl matrix is obtained in the preceding
step by the leaf peeling method. To this last problem the method from Section
\ref{Sect Local Problem} is applied, with $m_{1}$ being the number of the
edges of the star shaped graph. Obviously, in this last step there is no need
to compute functions (\ref{phiS}), so that the algorithm stops after
recovering the potential $q(x)$ on all the edges.

In the next section we discuss the numerical implementation of the method and
some numerical examples.

\section{Numerical examples\label{Sect Numerical}}

\begin{figure}
[ptb]
\begin{center}
\includegraphics[
height=2.5482in,
width=4.5261in
]%
{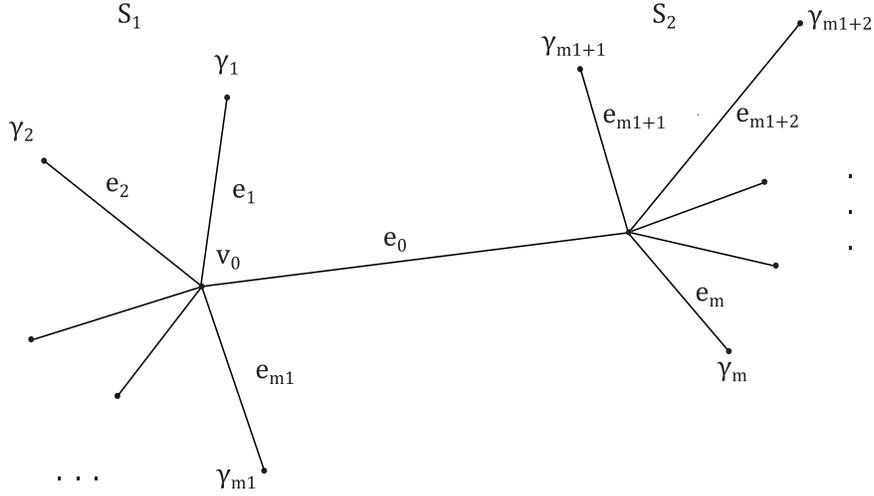}%
\caption{The tree graph considered in numerical tests.}%
\label{FigGraphExampleNumeric}%
\end{center}
\end{figure}

For numerical tests we consider the tree graph of the type depicted in Fig.
\ref{FigGraphExampleNumeric}.%

\textbf{Example 1. }Consider a graph from Fig. \ref{FigGraphExampleNumeric}
consisting of nine edges, $m=8$, $m_{1}=5$, the lengths of the edges
$e_{0},\ldots,e_{8}$ are
\begin{equation}
L_{0}=1.4,\,L_{1}=\frac{e}{2},\,L_{2}=1,\,L_{3}=\frac{\pi}{2},\,L_{4}%
=\frac{\pi}{3},\,L_{5}=\frac{e^{2}}{4},\,L_{6}=1.1,\,L_{7}=1.2,\,L_{8}=1.
\label{L}%
\end{equation}
The corresponding nine components of the potential are defined as follows%
\[
\,q_{0}(x)=J_{0}(9x)+1,\,q_{1}(x)=\left\vert x-1\right\vert +1,\,q_{2}%
(x)=e^{-(x-\frac{1}{2})^{2}},\,q_{3}(x)=\sin\left(  8x\right)  +\frac{2\pi}%
{3},\,
\]%
\[
q_{4}(x)=\cos\left(  9x^{2}\right)  +2,\,q_{5}(x)=\frac{1}{x+0.1}%
,\,q_{6}(x)=\frac{1}{\left(  x+0.1\right)  ^{2}},\,q_{7}(x)=e^{x},
\]%
\[
\,q_{8}(x)=\left\{
\begin{tabular}
[c]{ll}%
$-35.2x^{2}+17.6x,$ & $0\leq x<0.25$\\
$35.2x^{2}-35.2x+8.8,$ & $0.25\leq x<0.75$\\
$-35.2x^{2}+52.8x-17.6,$ & $0.75\leq x\leq1,$%
\end{tabular}
\ \ \ \ \right.  .
\]

The potential $q_{8}(x)$ (from \cite{Brown et al 2003}, \cite{Rundell Sacks})
is referred to below as saddle potential. To compute the Weyl matrix (direct
problem) at a number of points $\rho_{k}^{2}$, we followed the approach from
\cite{AKK2023}. We took $180$ points $\rho_{k}$ chosen according to the rule
$\rho_{k}=10^{\alpha_{k}}+0.1i$ with $\alpha_{k}$ being distributed uniformly
on $\left[  0,2\right]  $. Such choice delivers a set of points which are more
densely distributed near $\rho=1+0.1i$ and more sparsely near $\rho=100+0.1i$.
In all our computations $N$ was chosen as $N=9$ (ten coefficients
$g_{i,n}(L_{i})$ and ten coefficients $s_{i,n}(L_{i})$).

\begin{figure}
[ptb]
\begin{center}
\includegraphics[
height=3.8502in,
width=5.127in
]%
{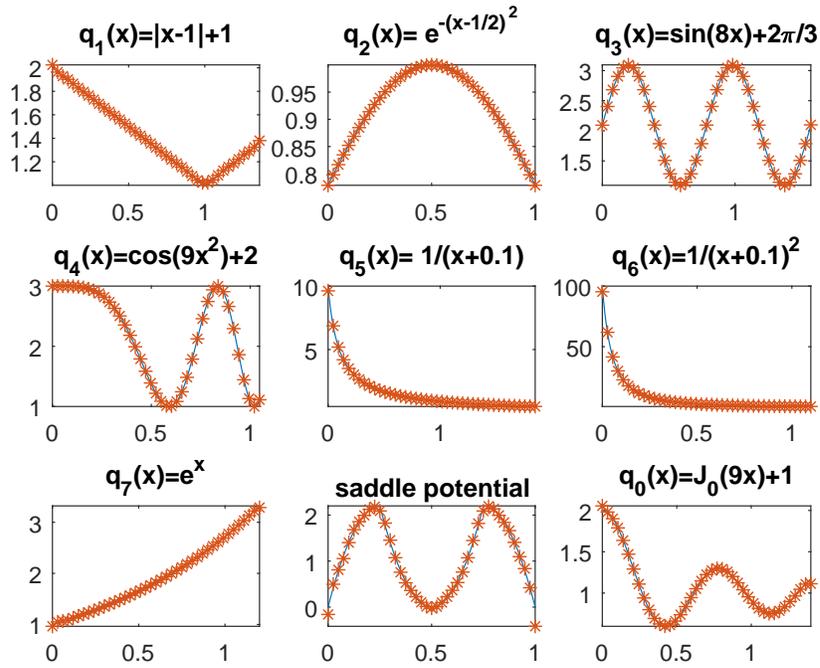}%
\caption{The potential of the quantum graph from Example 1, recovered from the
Weyl matrix given at 180 points, with $N=9$.}%
\label{FigGraph9}%
\end{center}
\end{figure}

In Fig. \ref{FigGraph9} we show the exact potentials (continuous line)
together with the recovered ones (marked with asterisks), computed with the
proposed algorithm.

\medskip%

Here the maximum relative error was attained in the case of the saddle
potential, and it resulted in approximately $0.085$ in the vicinity of the
endpoint $x=1$. All other potentials were computed considerably more accurately.

It is interesting to track how accurately the two spectra were computed on
each edge. For example, Table 1 presents some of the \textquotedblleft
exact\textquotedblright\ Dirichlet-Dirichlet eigenvalues on $e_{2}$ computed
with the aid of the Matslise package \cite{Ledoux et al} (first column), the
approximate eigenvalues, computed as described in subsection
\ref{Subsect reduction to two spectra} by calculating zeros of $S_{2,9}%
(\rho,L_{2})$ and the absolute error of each presented eigenvalue. The zeros
of $S_{2,9}(\rho,L_{2})$ were computed by converting this function into a
spline on the interval $\left(  0,650\right)  $ and using the Matlab command
\emph{fnzeros. }Notice that both the absolute and relative errors remain
remarkably small even for large indices.

\medskip

\bigskip%
\begin{tabular}
[c]{|l|l|l|l|}\hline
\multicolumn{4}{|l|}{Table 1: Dirichlet-Dirichlet eigenvalues of $q_{2}(x)$%
}\\\hline
$n$ & $\lambda_{n}$ & $\widetilde{\lambda}_{n}$ & $\left\vert \lambda
_{n}-\widetilde{\lambda}_{n}\right\vert $\\\hline
$1$ & $10.8381543818$ & $10.8381543825$ & $7.1\cdot10^{-10}$\\\hline
$11$ & $1195.1450218516$ & $1195.1450218560$ & $4.4\cdot10^{-9}$\\\hline
$51$ & $25671.7636244$ & $25671.7636252$ & $8.0\cdot10^{-7}$\\\hline
$101$ & $100680.7570614$ & $100680.7570623$ & $9.4\cdot10^{-7}$\\\hline
$201$ & $398742.8099714$ & $398742.8099723$ & $9.9\cdot10^{-7}$\\\hline
\end{tabular}

\bigskip\bigskip

Similar results were obtained for the Neumann-Dirichlet eigenvalues.

It is interesting to observe how the error caused by the leaf peeling
procedure is accumulated. For this we consider the second example, in which
the potentials from the first sheaf are repeated on the second one.

\begin{figure}
[ptb]
\begin{center}
\includegraphics[
height=3.3087in,
width=6.9461in
]%
{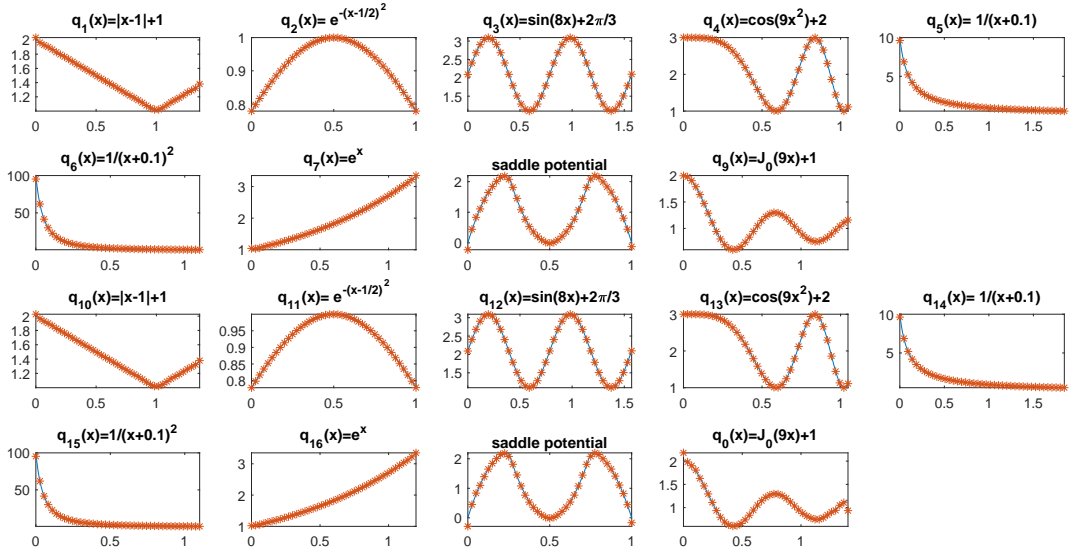}%
\caption{The potential of the quantum graph from Example 2, recovered from the
Weyl matrix given at 180 points, with $N=9$.}%
\label{FigGraph18}%
\end{center}
\end{figure}

\textbf{Example 2. }Consider a graph from Fig. \ref{FigGraphExampleNumeric}
consisting of eighteen edges, $m=17$, $m_{1}=9$, the lengths of the edges
$e_{1},\ldots,e_{9}$ are%
\[
L_{1}=\frac{e}{2},\,L_{2}=1,\,L_{3}=\frac{\pi}{2},\,L_{4}=\frac{\pi}%
{3},\,L_{5}=\frac{e^{2}}{4},\,L_{6}=1.1,\,L_{7}=1.2,\,L_{8}=1,\,L_{9}=1.4.
\]
The corresponding eight components of the potential $q_{1}(x),\ldots,q_{8}(x)$
are defined as in the previous example, and additionally, $q_{9}%
(x)=J_{0}(9x)+1$. Next, $q_{i+9}(x)=q_{i}(x)$ and $L_{i+9}=L_{i}$ for
$i=1,\ldots,8$, and $q_{0}(x)=q_{9}(x)$, $L_{0}=L_{9}$. That is, the nine
components of the potential on the leaf edges of the first sheaf are all the
potentials from the first example, the component of the potential on $e_{0}$
is that from Example 1, and eight components of the potential on the second
sheaf are $q_{1}(x),\ldots,q_{8}(x)$ from Example 1. All other parameters were
chosen as in Example 1. In Fig. \ref{FigGraph18} the result of the computation
is presented. It is noteworthy that there is no considerable difference in the
accuracy of the potentials recovered before applying the leaf peeling
procedure (the first nine potentials) and after (the other nine ones). The
largest difference is observed in the case of the potentials $q_{9}(x)$ and
$q_{0}(x)$. Indeed, the maximum relative error for $q_{9}(x)$ resulted in
approximately $0.003$ while for $q_{0}(x)$ in $0.175$ at the right endpoint.
Interestingly enough, the accuracy of the recovery of all other potentials
practically did not change after the leaf peeling procedure. For example, the
maximum relative error for $q_{5}(x)$ was $0.037583$ while for its
\textquotedblleft twin\textquotedblright\ $q_{14}(x)$ it resulted in
$0.037590$.%

Thus, the leaf peeling procedure does not lead to a considerable error
accumulation and clearly can be applied several times that allows one to solve
inverse problems on quite complicated tree graphs.

It is worth noting that the whole computation takes few seconds performed in
Matlab 2017 on a Laptop equipped with a Core i7 Intel processor.

\section{Conclusions\label{Sect Concl}}

A new method for solving inverse problems on quantum tree graphs, consisting
in the recovery of a potential from a Weyl matrix is developed. It is based on
the leaf peeling method and Neumann series of Bessel functions representations
for solutions of Sturm-Liouville equations. The given data are used for
solving local inverse problems on sheaves, which in turn are reduced to
separate two-spectra inverse Sturm-Liouville problems on leaf edges. The leaf
peeling method allows one to remove the edges where the potential is already
recovered and compute the Weyl matrix for the smaller tree. This combination
of the leaf peeling method with the approach based on the Neumann series of
Bessel functions representations results in a simple, direct and accurate
numerical algorithm. Its performance is illustrated by numerical examples.

\textbf{Funding }The research of Sergei Avdonin was supported in part by the
National Science Foundation, grant DMS 1909869, and by Moscow Center for
Fundamental and Applied Mathematics. The research of Vladislav Kravchenko was
supported by CONACYT, Mexico via the project 284470 and partially performed at
the Regional mathematical center of the Southern Federal University with the
support of the Ministry of Science and Higher Education of Russia, agreement 075-02-2022-893.

\textbf{Data availability} The data that support the findings of this study
are available upon reasonable request.

\textbf{Declarations}

\textbf{Conflict of interest} The authors declare no competing interests.


\begin{thebibliography}{99}                                                                                               %


\bibitem {AbramowitzStegunSpF}Abramovitz M. and Stegun I. A. (1972),
\textit{Handbook of mathematical functions}, New York: Dover.

\bibitem {ArioliBenzi2018}Arioli M. and Benzi M. (2018), A finite element
method for quantum graphs. \textit{IMA J. Numer. Anal.}, \textbf{38}, no. 3, 1119--1163.

\bibitem {AvdoninBelinskiyMatthews}{Avdonin S., Belinskiy B. and Matthews J.}
(2011), {Inverse problem on the semi-axis: local approach}, \textit{Tamkang
Journal of Mathematics}, \textbf{42}, no. 3, 1--19.

\bibitem {AvdoninBell2015}Avdonin S. and Bell J. (2015), Determining physical
parameters for a neuronal cable model defined on a tree graph, \textit{Journal
of Inverse Problems and Imaging,} \textbf{9}, no. 3, 645-659.

\bibitem {AChoque}{Avdonin }S., Choque Rivero A., Leugering G. and Mikhaylov
V. (2015), On the inverse problem of the two velocity tree-like graph,
\emph{Zeit. Angew. Math. Mech.,} \textbf{95} , no. 12, 1490--1500.

\bibitem {AKK2023}Avdonin S. A., Khmelnytskaya K. V. and Kravchenko V. V.
Recovery of a potential on a quantum star graph from Weyl's matrix. arXiv:2210.15536.

\bibitem {AvdoninKravchenko2022}Avdonin S. A. and Kravchenko V. V. (2023),
Method for solving inverse spectral problems on quantum star graphs.
\emph{Journal of Inverse and Ill-posed Problems}, \textbf{31}, no. 1, 31-42.

\bibitem {AvdoninKurasov2008}Avdonin S. and Kurasov P. (2008), {Inverse
problems for quantum trees,} \textit{Inverse Problems and Imaging},
\textbf{2}, no. 1, 1--21.

\bibitem {AvdoninLeugeringMikhaylov2010}Avdonin S., Leugering G. and Mikhaylov
V. (2010), {On an inverse problem for tree-like networks of elastic strings},
\textit{Zeit. Angew. Math. Mech.}, \textbf{90}, no. 2, 136--150.

\bibitem {AvdoninZhao2020}{Avdonin} S. and Zhao Yu. (2021), Leaf peeling
method for the wave equation on metric tree graphs. \emph{Inverse Probl.
Imaging} \textbf{15}, no. 2, 185--199.

\bibitem {Baricz et al Book}Baricz A., Jankov D. and Pog\'{a}ny T. K. (2017),
\emph{Series of Bessel and Kummer-type functions. }Lecture Notes in
Mathematics, 2207. Springer, Cham.

\bibitem {Belishev Vakulenko 2006}Belishev M. and Vakulenko A. (2006), Inverse
problems on graphs: Recovering the tree of strings by the BC-method, \emph{J.
Inv. Ill-Posed Problems}, \textbf{14} , 29-46.

\bibitem {BerkolaikoKuchment}Berkolaiko G. and Kuchment P. (2013),
\textit{Introduction to Quantum Graphs}, AMS, Providence, R.I.

\bibitem {Brown et al 2003}Brown B. M., Samko V. S., Knowles I. W., Marletta
M. (2003), Inverse spectral problem for the Sturm--Liouville equation,
\emph{Inverse Probl.} \textbf{19}, 235--252.

\bibitem {Chadan et al 1997}Chadan Kh., Colton D., P\"{a}iv\"{a}rinta L.,
Rundell W. (1997), \emph{An introduction to inverse scattering and inverse
spectral problems}. SIAM, Philadelphia.

\bibitem {KarapetyantsKravchenkoBook}Karapetyants A. N. and Kravchenko V. V.
(2022), \emph{Methods of mathematical physics: classical and modern.}
Birkh\"{a}user, Cham.

\bibitem {Kr2019JIIP}Kravchenko V. V. (2019), On a method for solving the
inverse Sturm--Liouville problem,\textbf{ }\emph{J. Inverse Ill-posed Probl.
}\textbf{27}, 401--407.

\bibitem {KrBook2020}Kravchenko V. V. (2020), \emph{Direct and inverse
Sturm-Liouville problems: A method of solution}, Birkh\"{a}user, Cham.

\bibitem {Kr2022Completion}Kravchenko V. V. (2022), Spectrum completion and
inverse Sturm-Liouville problems. Math Meth Appl Sci.; 1-15. doi:10.1002/mma.8869.

\bibitem {KKC2022Mathematics}Kravchenko V. V., Khmelnytskaya K. V. and
\c{C}etinkaya F. A. (2022), Recovery of inhomogeneity from output boundary
data, \emph{Mathematics}, \textbf{10}, 4349, https://doi.org/10.3390/math10224349.

\bibitem {KNT}Kravchenko V. V., Navarro L. J. and Torba\ S. M. (2017),
Representation of solutions to the one-dimensional Schr\"{o}dinger equation in
terms of Neumann series of Bessel functions, \emph{Appl. Math. Comput.}
\textbf{314}, 173--192.

\bibitem {KT2015JCAM}Kravchenko V. V. and Torba\ S. M. (2015), Analytic
approximation of transmutation operators and applications to highly accurate
solution of spectral problems,\emph{ Journal of Computational and Applied
Mathematics }\textbf{275}, 1-26.

\bibitem {KT2021 IP1}Kravchenko V. V. and Torba S. M. (2021), A direct method
for solving inverse Sturm-Liouville problems, \emph{Inverse Probl.
}\textbf{37}, 015015 (32pp).

\bibitem {KT2021 IP2}Kravchenko V. V. and Torba S. M. (2021), A practical
method for recovering Sturm-Liouville problems from the Weyl function,
\emph{Inverse Probl. }\textbf{37}, 065011 (26pp).

\bibitem {Kurasov et al 2005}Kurasov P. and Nowaczyk M. (2005), Inverse
spectral problem for quantum graphs, \emph{J. Phys. A.,} \textbf{38}, 4901-4915.

\bibitem {Ledoux et al}Ledoux V., Daele M.V. and Berghe G.V. \ (2005),
MATSLISE: a MATLAB package for the numerical solution of Sturm--Liouville and
Schr\"{o}dinger equations, \emph{ACM Trans. Math. Softw.} \textbf{31}, 532--554.

\bibitem {LevitanInverse}Levitan B. M. (1987), \emph{Inverse Sturm-Liouville
problems}, VSP, Zeist.

\bibitem {Marchenko}Marchenko V. A. (2011), \emph{Sturm-Liouville operators
and applications: revised edition}, AMS Chelsea Publishing.

\bibitem {Mugnolo}Mugnolo D. (2014), \textit{Semigroup Methods for Evolution
Equations on Networks,} Understanding Complex Systems, Springer, Cham.

\bibitem {Rundell Sacks}Rundell W. and Sacks P. E. (1992), Reconstruction
techniques for classical inverse Sturm--Liouville problems, \emph{Math.
Comput. }\textbf{58}, 161--183.

\bibitem {SitnikShishkina Elsevier}Shishkina E. L. and Sitnik S. M. (2020),
\emph{Transmutations, singular and fractional differential equations with
applications to mathematical physics}, Elsevier, Amsterdam.

\bibitem {SavchukShkalikov}Savchuk A. M., Shkalikov A. A. (2005), Inverse
problem for Sturm--Liouville operators with distribution potentials:
reconstruction from two spectra, \emph{Russ. J. Math. Phys. }\textbf{12}, 507--514.

\bibitem {Watson}Watson G. N. (1996), \emph{A Treatise on the theory of Bessel
functions, 2nd ed., reprinted}, Cambridge University Press, Cambridge.

\bibitem {Wilkins}Wilkins J. E. (1948), Neumann series of Bessel functions.
\emph{Trans. Amer. Math. Soc.} \textbf{64}, 359--385.

\bibitem {Yurko2005}Yurko V. A. (2005), Inverse Sturm-Lioville operator on
graphs, \emph{Inverse Problems}, \textbf{21}, 1075-1086.

\bibitem {Yurko2007}Yurko V. A. (2007), \emph{Introduction to the theory of
inverse spectral problems}, Fizmatlit, Moscow, (in Russian).
\end{thebibliography}
\end{document}